\documentclass{amsart}
\usepackage{amssymb,amscd,amsthm,amsmath}
\usepackage{framed}
\usepackage{ascmac}
\usepackage{enumerate}
\usepackage{graphicx}
\usepackage{setspace}
\usepackage{bm,bbm}
\usepackage{booktabs}

\theoremstyle{plain}
\newtheorem{theorem}{Theorem}[section]

\newtheorem{corollary}[theorem]{Corollary}

\theoremstyle{definition}
\newtheorem{definition}[theorem]{Definition}
\newtheorem{remark}[theorem]{Remark}
\newtheorem{example}[theorem]{Example}

\newtheorem*{acknowledgment}{Acknowledgment}

\theoremstyle{remark}

\numberwithin{equation}{section}

\newcommand{\Z}{\mathbbm{Z}}
\newcommand{\Q}{\mathbbm{Q}}

\newcommand{\Spl}[1]{\mathrm{Spl}_\Q\left(#1\right)}
\newcommand{\Gal}{\mathrm{Gal}}
\newcommand{\Ker}{\mathrm{Ker}\,}
\newcommand{\Hom}{\mathrm{Hom}_{\mathrm{cont}}}

\newcommand{\Est}{E_{t,s}}
\newcommand{\Estast}{E_{t,s}^{\ast}}

\newcommand{\Ep}{\mathcal{E}_p}
\newcommand{\Epast}{\mathcal{E}_p^{\ast}}
\newcommand{\Epr}{\mathcal{E}_{p,r}}
\newcommand{\Eprast}{\mathcal{E}_{p,r}^{\ast}}
\newcommand{\Bru}{\mathrm{Bru}}
\newcommand{\Lec}{\mathrm{Lec}}

\newcommand{\Tor}{\mathrm{tor}}

\title{On Lecacheux's family of quintic polynomials}

\author[A. Hoshi]{Akinari Hoshi}
\address{Department of Mathematics, Niigata University, Niigata 950-2181, Japan}
\email{hoshi@math.sc.niigata-u.ac.jp}

\author[M. Koshiba]{Masakazu Koshiba}
\address{Graduate School of Science and Technology, Niigata University, Niigata 950-2181, Japan}
\email{koshiba@m.sc.niigata-u.ac.jp}

\thanks{
This work was partially supported by JSPS KAKENHI Grant Number 19K03418. 
}


\begin{document}
\begin{abstract}
Kida, Rikuna and Sato \cite{KRS10} developed a classification theory 
for Brumer's quintic polynomials via Kummer theory arising 
from associated elliptic curves. 
We generalize their results to elliptic curves associated to 
Lecacheux's quintic $F_{20}$-polynomials instead of Brumer's quintic 
$D_5$-polynomials. 
\end{abstract}

\maketitle
\tableofcontents

\section{Introduction}\label{seInt}
Let $K$ be a field with ${\rm char}\, K\neq 2,\ 5$ 
and $C_n$ be the cyclic group of order $n$. 
Let $D_5\simeq C_5\rtimes C_2$ be the dihedral group of order $10$ and 
$F_{20}\simeq C_5\rtimes C_4$ be the Frobenius group of order $20$. 
Let $K(s,t)$ be the rational function field over $K$ with two variables $s,t$. 
Brumer's quintic polynomial $\Bru(t,s;X)$ is defined to be 
\begin{align}
\Bru(t,s;X) := X^5&+ (t-3)X^4 - (t - s - 3)X^3\label{Bru}\\
&+ (t^2 - t - 2s - 1)X^2 + sX + t \in K(s,t)[X].\nonumber
\end{align}

The polynomial $\Bru(t,s;X)$ is $K$-generic for $D_5$, 
namely 
{\rm (i)} the Galois group of $\Bru(t,s;X)$ over $K(s,t)$ 
is isomorphic to $D_5$; and 
{\rm (ii)} every $D_5$-Galois extension $L/M$, $\# M=\infty$, $M\supset K$, 
can be obtained as $L=\mathrm{Spl}_M (\Bru(b,a;X))$, 
the splitting field of $\Bru(b,a;X)$ over $M$, for some $a,b\in M$ 
(see Jensen, Ledet and Yui \cite[Theorem 2.3.5]{JLY02}). 

Kida, Rikuna and Sato \cite{KRS10} investigated 
Brumer's quintic $\Bru(t,s;X)$ via Kummer theory 
arising from elliptic curves. 
The splitting field $\mathrm{Spl}_{K(s,t)}{(\Bru(t,s;X))}$ 
of $\Bru(t,s;X)$ over $K(s,t)$ contains the unique quadratic 
subfield $K(s,t)(\sqrt{d_{t,s}})$ where 
\begin{align}
d_{t,s} :=&-4s^3+(t^2-30t+1)s^2+2t(3t+1)(4t-7)s\label{dts}\\
&-t(4t^4-4t^3-40t^2+91t-4)\in K(s,t).\nonumber
\end{align}

In this paper, we study the case where $K=\Q$. 
We search elements $\alpha$ and $\beta$ in $\Q(s,t)$ such that 
the quadratic subfields of 
${\rm Spl}_{\Q(s,t)}({\Bru(\beta,\alpha;X)})$ 
and of ${\rm Spl}_{\Q(s,t)}(\Bru(t,s;X))$ coincide. 
According to Kida, Rikuna and Sato \cite[Section 2]{KRS10}, 
we restrict ourselves to treat the case $\beta=t$ 
and consider the equation 
\begin{align*}
d_{t,s}u^2=d_{t,\alpha}.
\end{align*} 
Define 
\begin{align*}
d=d_{t,s},\ x=-4d\alpha,\ y=4d^2u.
\end{align*}
Then we obtain the associated elliptic curve 
\begin{align}
\Est : y^2 =x^3 &+ d(t^2-30t+1)x^2 - 8d^2t(3t+1)(4t-7)x\label{est}\\
&-16d^3t(4t^4-4t^3-40t^2+91t-4)\nonumber
\end{align}
to Brumer's quintic polynomial $\Bru(t,s;X)$. 
This elliptic curve $\Est$ has an isogeny $\phi$ of degree $5$ 
defined over $\Q(s,t)$. 
The $5$-division polynomial of $\Est$ 
(see Silverman \cite[Exercise 3.7]{Sil86}) 
has a quadratic factor 
\begin{align*}
f_2(x)=x^2+s_1x+\frac{s_2}{5}
\end{align*} 
where 
\begin{align*}
s_1&=-16t^7+96t+(96s+64)t^5+(4s^2-616s-1148)t^4+(-140s^2+720s+1996)t^3\\
&\quad \ +(-16s^3+608s^2+144s-444)t^2+(80s^3-140s^2-56)t-16s^3+4s^2,\\
s_2&=256t^{14}-230t^{13}+(-3072s-256)t^{12}+(-128s^2+28928s+44160)t^{11}\\
&\quad \ +(14080s^2-36096s-113920)t^{10}+(1280s^3-121600s^2-227328s-55936)t^9\\
&\quad \ +(16s^4-33600s^3+223200s^2+866496s+696464)t^8\\
&\quad \ +(-4144s^4+205440s^3-14560s^2-1512960s-1663984)t^7\\
&\quad \ +(-128s^5+47072s^4-292480s^3-333792s^2+2181696s+2733120)t^6\\
&\quad \ +(4736s^5-147072s^4+387200s^3+1014560s^2-1471872s-2641440)t^5\\
&\quad \ +(256s^6-27520s^5+141920s^4-423040s^3-1920800s^2-578496s+348304)t^4\\
&\quad \ +(-1792s^6+18560s^5-338928s^4-408000s^3+169600s^2+68672s-16256)t^3\\
&\quad \ +(-4608s^6+6144s^5+544s^4-960s^3+128s^2)t+256s^6-128s^5+16s^4.
\end{align*}
Take a root $\theta$ of $f_2(x)=0$. 
Then we obtain a point $A\in\Est(\overline{\Q(s,t)})$ of order $5$ 
with $x(A)=\theta$. 
Apply the V\'{e}lu formula \cite{Vel71} to $\langle A \rangle$ 
and take $\Estast = \Est/\langle A \rangle$ as the image of $\phi$ 
(see Kida, Rikuna and Sato \cite[Section 2]{KRS10}): 
\begin{align}
\Estast : y^2 = x^3&+d(t^2-30t+1)x^2-8d^2(26t^4-310t^3+327t^2+315t+26)x
\label{Eastts}\\
&+16d^3(68t^6-1120t^5+3804t^4+1760t^3+6929t^2+1380t+68).\nonumber
\end{align}

After the specialization $\Q(s,t)^2\ni (s,t)\mapsto 
(s^\prime,t^\prime)\in \Q^2$, 
we obtain that 
$\Bru(t^\prime,s^\prime;X)$, $E_{t^\prime,s^\prime}$ and 
$E^\ast_{t^\prime,s^\prime}$ 
are defined over $\Q$. 
After the specialization, for $s,t\in \Q$, 
we also write $\Bru(t,s;X)$, $\Est$ and $\Estast$
which are defined over $\Q$ (not $\Q(s,t)$). 

Let
\begin{align*}
\phi^{\ast} : \Estast \rightarrow \Est
\end{align*}
be the dual isogeny of $\phi$. 
Then the quotient group $\Est(\Q)/\phi^{\ast}(\Estast(\Q))$ 
is finite 
by weak Mordel--Weil theorem (see \cite[Chapter VIII, \S1]{Sil86}). 

\begin{definition}[{Kida, Rikuna and Sato \cite[page 694]{KRS10}}]
\label{defBru}
Let $s,t$ be rational numbers. 
For each $\Q$-rational point $P = (x(P),y(P)) \in \Est(\Q)$, 
Brumer's polynomial $\Bru(P;X)$ with respect to $P$ is defined to be 
\begin{align*}
\Bru(P;X) := \Bru \left(t,\frac{x(P)}{-4d};X\right)
\end{align*}
where $\Bru(t,s;X)$ is Brumer's polynomial as in (\ref{Bru}) 
and $d=d_{t,s}$ is given as in (\ref{dts}).
\end{definition}
We remark that 
there exists 
a rational point 
$P_0=(-4ds,4d^2)\in \Est(\Q)$ 
and 
by the definition
we have $\Bru(P_0;X)=\Bru(t,s;X)$.
\begin{theorem}[{Kida, Rikuna and Sato \cite[Theorem 2.1]{KRS10}}]
\label{thKRS}
Let $s,t$ be rational numbers. 
Let $\Est$ be the elliptic curve as in $(\ref{est})$. 
Let $\Bru(P;X)$ be Brumer's polynomial with respect to $P$ 
as in Definition \ref{defBru} 
with the splitting field $\Spl{\Bru(P;X)}$ over $\Q$.\\
{\rm (i)} 
For any $\Q$-rational point $P \in \Est (\Q)$, 
$\Bru(P;X)$ is reducible over $\Q$ 
if and only if $P \in \phi^{\ast}(\Estast (\Q))$;\\
{\rm (ii)} 
There exists a bijection between the following two finite sets
\begin{center}
\{subgroup of order $5$ in $\Est (\Q)/\phi^{\ast}(\Estast(\Q))$\}
\end{center}
and
\begin{center}
\{$\Spl{\Bru(P;X)}\mid P \in \Est (\Q)\setminus \phi^{\ast}(\Estast(\Q))$\}.
\end{center}

The bijection is induced by the correspondence 
$\Est (\Q)\ni P \mapsto \Spl{\Bru(P;X)}$.
\end{theorem}

The aim of this paper is to generalize Theorem \ref{thKRS} 
to elliptic curves associated to Lecacheux's quintic $F_{20}$-polynomial 
$\Lec(p,r;X)$ instead of Brumer's quintic $D_5$-polynomial $\Bru(t,s;X)$. 

Let $\Q(p,r)$ be the rational function field over $\Q$ with 
two variables $p,r$. 
Lecacheux's quintic polynomial $\Lec(p,r;X)$ is defined to be 
\begin{align}
\Lec(p,r;X) := &\ X^5+\left(r^2(p^2+4)-2p-\dfrac{17}{4}\right) X^4
+\left(3r(p^2+4)+p^2+\dfrac{13}{2}p+5\right) X^3\label{Lec}\\
&-\left(r(p^2+4)+\dfrac{11}{2}p-8\right)X^2+(p-6)X+1\in\Q(p,r)[X].\nonumber
\end{align} 
The polynomial $\Lec(p,r;X)$ is known to be $\Q$-generic for $F_{20}$ 
(see Jensen, Ledet and Yui \cite[Theorem 2.3.6]{JLY02}). 

We will define the elliptic curve $\Epr$ 
associated to Lecacheux's quintic polynomial $\Lec(p,r;X)$. 
Define 
\begin{align}
W_{p,r} &:= 16(p^2+4)r^3 +4(p^2+4)r^2-4(19p+41)r-16p-199,\label{Wpr}\\
D_{p,r} &:= W_{p,r}\left((p^4+5p^2+4)+p(p^2+3)\sqrt{p^2+4}\right)
\Big/8.\nonumber
\end{align}

The splitting field $\mathrm{Spl}_{\Q(p,r)}(\Lec(p,r;X))$ 
contains the unique quadratic (resp. quartic) subfield 
$\Q(p,r)(\sqrt{p^2+4})$ (resp. $\Q(p,r)(\sqrt{D_{p,r}})$) 
(see Hoshi and Miyake \cite[Lemma 7.3 and Lemma 7.4]{HM10}; 
$\Lec(p,r;X)$ is $g^{F_{20}}_{p,r}(X)$ in \cite{HM10}).

We search $\beta$ such that the quartic subfields of $\Spl{\Lec(p,\beta;X)}$ 
and of $\Spl{\Lec(p,r;X)}$ coincide. 
We consider the equation 
\begin{align*}
D_{p,r}u^2=D_{p,\beta}.
\end{align*}
Write $D= D_{p,r}$ and $W=W_{p,r}$. 
Then the above equation becomes 
\begin{align*}
Wu^2=W_{p,\beta}.
\end{align*}
Define
\begin{align*}
x:=4(p^2+4)W\beta,\ y:=2(p^2+4)W^2u. 
\end{align*}
Then we get the associated elliptic curve
\begin{align}
\Epr:y^2&=x^3+(p^2+4)Wx^2-4(19p+41)(p^2+4)W^2x-4(p^2+4)^2(16p+199)W^3
\label{epr}
\end{align}
to Lecacheux's quintic polynomial $\Lec(p,r;X)$. 

The curve $\Epr$ has an isogeny $\nu$ of degree $5$ defined over $\Q(p,r)$. 
We see that the $5$-division polynomial of $\Epr$ 
(see Silverman \cite[Exercise 3.7]{Sil86}) 
has the quadratic factor 
\begin{align*}
f_2(x)=x^2+s_1x+\frac{s_2}{5}
\end{align*}
where 
\begin{align*}
s_1&=(64r^3+16r^2)p^4+(-304r-64)p^3+(512r^3+128r^2-656r-796)p^2\\
&\quad \ +(-1216r-256)p+1024r^3+256r^2-2624r-3184,\\
s_2&=(4096r^6+2048r^5+256r^4)p^8\\
&\quad \ +(22528r^6+ 11264r^5 -37504r^4 -17920r^3 -2048r^2)p^7\\
&\quad \ +( -54272r^6 -27136r^5 -301376r^4 -221440r^3+ 55680r^2+ 38912r+4096)p^6\\
&\quad \ +(270336r^6+ 135168r^5+ 226304r^4 -366720r^3+ 802368r^2+ 781952r+124416)p^5\\
&\quad \ +(-1044480r^6-522240r^5-328960r^4+1331200r^3+ 1009440r^2+ 3341120r+1106960)p^4\\
&\quad \ +( 1081344r^6+540672r^5+ 3610624r^4 -2073600r^3 -4815360r^2-5468640r+1409880)p^3\\
&\quad \ +( -4702208r^6-2351104r^5+ 11899904r^4+ 21278720r^3 -5146432r^2 -13628000r- 11636500)p^2\\
&\quad \ +(1441792r^6+ 720896r^5+ 8421376r^4 -3573760r^3 -32230400r^2-34385792r+3648864)p\\
&\quad \ -6619136r^6 -3309568r^5+ 33509376r^4+ 49643520r^3 -33173248r^2 -105479552r- 63995216.
\end{align*}
Take a root $\theta$ of $f_2(x)=0$. 
Then we obtain a point $A\in\Epr(\overline{\Q(p,r)})$ of order $5$ 
with $x(A)=\theta$, 
$\Eprast = \Epr/\langle A \rangle$ as the image of $\nu$ 
and the dual isogeny $\nu^{\ast} : \Eprast \rightarrow \Epr$ of $\nu$ 
as in (\ref{Eastts}) (see also Kida, Rikuna and Sato \cite[Section 2]{KRS10}): 
\begin{align*}
\Eprast : y^2 = x^3&+(p^2+4)W x^2-4(p^2+4)(52p^2-625p+833)W^2x\\
& +4(p^2+4)^2(272p^2-5000p+21713) W^3.
\end{align*}

As in the case of Brumer's quintic, 
after the specialization $\Q(p,r)^2\ni (p,r)\mapsto 
(p,r)\in \Q^2$, we also write 
$\Lec(p,r;X)$, $\Epr$ and $\Eprast$ for $p,r\in \Q$ 
which are defined over $\Q$ (not $\Q(p,r)$). 

%
\begin{definition}\label{defLec}
Let $p,r$ be rational numbers. 
For each $\Q$-rational point $P=(x(P),y(P)) \in \Epr(\Q)$, 
Lecacheux's polynomial $\Lec(P;X)$ with respect to $P$ is defined to be 
\begin{align*}
\Lec(P;X) := \Lec\left(p,\frac{x(P)}{4(p^2+4)W};X\right)
\end{align*}
where $\Lec(p,r;X)$ is Lecacheux's polynomial as in (\ref{Lec}) 
and $W=W_{p,r}$ is given as in (\ref{Wpr}). 
\end{definition}
We note that 
there exists the point $Q_0=(4r(p^2+4)W,2(p^2+4)W^2)\in \Epr(\Q)$ 
and we have $\Lec(Q_0;X)=\Lec(p,r;X)$ by the definition. 

The main theorem of this paper can be described as follows: 
\begin{theorem}\label{t1.4}
Let $p,r$ be rational numbers. 
Let $\Epr$ be the elliptic curve as in $(\ref{epr})$. 
Let $\Lec(P;X)$ be Lecacheux's polynomial with respect to $P$ 
as in Definition \ref{defLec} with the splitting field 
$\Spl{\Lec(P;X)}$ over $\Q$.\\
{\rm (i)} For any $\Q$-rational point $P \in \Epr (\Q)$, 
$\Lec(P;X)$ is reducible over $\Q$ 
if and only if $P \in \nu^{\ast}(\Eprast (\Q))$;\\
{\rm (ii)} 
There exists a bijection between the following two finite sets
\begin{center}
\{subgroup of order $5$ in $\Epr (\Q)/\nu^{\ast}(\Eprast(\Q))$\}
\end{center}
and
\begin{center}
\{$\Spl{\Lec(P;X)}\mid P \in \Epr (\Q)\setminus \nu^{\ast}(\Eprast(\Q))$\}.
\end{center}

The bijection is induced by the correspondence 
$\Epr (\Q)\ni P \mapsto \Spl{\Lec(P;X)}$. 
\end{theorem}

\section{Constructions of $\Bru(t,s;X)$ and $\Lec(p,r;X)$}
We recall constructions of Brumer's and Lecacheux's polynomials 
in Lecacheux \cite[pages 209--214]{Lec98}. 

\subsection{Construction of Brumer's polynomial $\Bru(t,s;X)$} 
We consider the elliptic curve:
\begin{align*}
E_t^{\ast}: y^2+(1-t)xy-ty=x^3-tx^2
\end{align*}
with $5$-torsion points
\begin{align*}
A=(0,0),\ 2A=(t,t^2),\ 3A=(t,0),\ 4A=(0,t).
\end{align*}
The curve $E_t^{\ast}$ is also called Tate normal form 
(see Husem\"oller \cite[page 93, Definition 4.1]{Hus04}). 
The $j$-invariant of $E_t^{\ast}$ is 
\begin{align*}
\frac{(t^4-12t^3+14t^2+12t+1)^3}{t^5(t^2-11t-1)}.
\end{align*}
There exists the elliptic curve 
$E_t=E_t^{\ast}/\langle A \rangle$ up to isomorphism 
with the isogeny 
\begin{align*}
\phi:&\ E_t^{\ast}\rightarrow E_t,\\
&\ X=\frac{t}{x} \mapsto X^\prime
=\frac{2(X-2)(X^2+2Xt-1)(2X^2-2Xt-2X+t)}{X(X-1)^2}
\end{align*}
of degree $5$. 
Then 
by solving this for $X$, we have 
\begin{align*}
X^5&+(t-3)X^4+(1-\frac{1}{4}X^\prime-2t^2-\frac{7}{2}t)X^3\\
&+(4t+3+5t^2+\frac{1}{2}X^\prime)X^2
+(-2t^2-2-\frac{1}{4}X^\prime-\frac{5}{2}t)X+t=0.
\end{align*}
Define 
\begin{align*}
s= -2t^2-2-\frac{1}{4}X^\prime-\frac{5}{2}t.
\end{align*}
Then the left-hand side of this equation becomes
\begin{align*}
\Bru(t,s;x)=x^5+(t-3)x^4-(t-s-3)x^3+(t^2-t-2s-1)x^2+sx+t.
\end{align*}

We find that the elliptic curve $E_t$ and 
the elliptic curve $\Est$ associated to $\Bru(t,s;X)$ 
as in $(\ref{est})$ are isomorphic 
over some extension field 
(see also Kida, Rikuna and Sato \cite[page 695]{KRS10}). 
The $j$-invariants of $E_t$ and of $\Est$ are the same 
\begin{align*}
\frac{(t^4+228t^3+494t^2-228t+1)^3}{t(t^2-11t-1)^5}.
\end{align*}
\subsection{Construction of Lecacheux's polynomial $\Lec(p,r;X)$} 
We consider the elliptic curve
\begin{align*}
\Epast:y^2-\frac{1}{4}(p^2+4)(x^2+1)=\frac{1}{2}(x^2-px-1)(2x-p).
\end{align*}
The elliptic curve $\Epast$ has $5$-torsion points:
\begin{align*}
A=(\alpha , \beta),\ 
2A=\left(-\frac{1}{\alpha} ,\ \frac{\beta}{\alpha} \right),\ 
3A=\left(-\frac{1}{\alpha} ,\ -\frac{\beta}{\alpha} \right),\ 
4A=(\alpha , -\beta)
\end{align*}
where 
$\alpha$ and $-1/\alpha$ are roots of $x^2-px-1$ and $\beta$ satisfies 
\begin{align*}
\beta^2 = \frac{1}{4}(p^2+4)(\alpha^2+4)=\frac{1}{4}(p^2+4)^{\frac{3}{2}}\alpha.
\end{align*}
The $j$-invariant of $\Epast$ is 
\begin{align*}
\frac{(p^2-12p+16)^3}{p-11}. 
\end{align*}
There exists the elliptic curve $\Ep=\Ep^{\ast}/\langle A \rangle$ 
up to isomorphism with the isogeny 
\begin{align*}
\phi:&\ \Ep^{\ast}\rightarrow \Ep\\
&\ x \mapsto r=\frac{x+2p}{p^2+4}+\frac{(p^2+4)(px+2)}{L^2}
+\frac{x(p+2)+(p^2-p+6)}{L}-\frac{5p}{2(p^2+4)}
\end{align*}
of degree $5$ where $L=x^2-px-1$. 
Define $l=-L/(p^2+4)$. 
Solving the equation for $l$, we have
\begin{align*}
l^5&+\left(r^2(p^2+4)-2p-\frac{17}{4}\right)l^4
+\left(3r(p^2+4)+p^2+\frac{13}{2}p+5\right)l^3\\
&-\left(r(p^2+4)+\frac{11}{2}p-8\right)l^2+(p-6)l+1=0.
\end{align*}
The left-hand side of this equation yields  
\begin{align*}
\Lec(p,r;l)=l^5&+\left(r^2(p^2+4)-2p-\dfrac{17}{4} \right) l^4
+\left( 3r(p^2+4)+p^2+\dfrac{13}{2}p+5 \right) l^3\\
&-\left( r(p^2+4)+\dfrac{11}{2}p-8 \right) l^2+(p-6)l+1.
\end{align*}
The elliptic curve $\Ep$ and 
the associated elliptic curve $\Epr$ to $\Lec(p,r;X)$ 
as in $(\ref{epr})$ are isomorphic 
over some extension field 
with the $j$-invariant 
\begin{align*}
\frac{(p^2 + 228p + 496)^3}{(p-11)^5}.
\end{align*}

\section{Proof of Theorem \ref{t1.4}}
The idea of the proof of Theorem \ref{t1.4} is to 
combine the results given in Hoshi and Miyake \cite{HM10} and 
Kida, Rikuna and Sato \cite{KRS10}. 
According to \cite[page 1078, Equation (25)]{HM10}, for $p,r\in \Q$, 
we define $k=\Q(\sqrt{p^2+4})$ and 
\begin{align}
s&=-\frac{1}{4}(5p+8r+2p^2r+(2pr+5)\sqrt{p^2+4}),\label{eqHM}\\
t&=\frac{1}{2}(p+\sqrt{p^2+4}).\nonumber
\end{align}
Then it follows that 
${\rm Spl}_k(\Bru(t,s;X))={\rm Spl}_\Q(\Lec(p,r;X))$. 
The associated elliptic curves 
$\Est$ and $\Estast$ given as in (\ref{est}) and (\ref{Eastts}) 
are defined over $k$. 
According to \cite[Section 3]{KRS10}, 
we take elliptic curves $E_t$ and $E_t^{\ast}$ defined over $k$ by
\begin{align*}
E_t &: y^2-(t-1)xy-ty=x^3-tx^2-5t(t^2+2t-1)x-t(t^4+10t^3-5t^2+15t-1),\\
E_t^{\ast} &: y^2-(t-1)xy -ty=x^3-tx^2.
\end{align*}
The curves $\Est$ (resp. $\Estast$) and $E_t$ (resp. $E_t^{\ast}$) 
are isomorphic over $F=k(\sqrt{d_{t,s}})$ 
where $d_{t,s}$ is given in $(\ref{dts})$ and 
we take an isogeny $\phi:\Est \rightarrow \Estast$ and 
the dual isogeny $\phi^{\ast}: \Estast \rightarrow \Est$. 
We also take an isogeny 
$\lambda^{\ast} : E_t^{\ast} \rightarrow E_t$ of degree $5$. 
By \cite[Theorem 3.1]{KRS10}, there exists an injective homomorphism 
\begin{align*}
\Est(k)/\phi^{\ast}(\Estast(k)) \hookrightarrow
\Hom(\Gal(\overline{F}/F),\Ker \lambda^{\ast}(k))).
\end{align*}

We will prove that there exists an injective homomorphism
\begin{align*}
\Epr(\Q)/\nu^{\ast}(\Eprast(\Q)) \hookrightarrow
\Hom(\Gal(\overline{F}/F),\Ker \lambda^{\ast}(k)).
\end{align*}
We see that 
the elliptic curves $\Epr$ and $\Est$ are isomorphic over $k$ 
with $j$-invariant 
\begin{align*}
\frac{(p^2+228p+496)^3}{(p-11)^5}.
\end{align*}
Indeed, we may find an isomorphism $f:\Epr \rightarrow \Est$ 
which is given explicitly as 
\begin{align*}
(x,y) \mapsto (ax+b,uy)
\end{align*}
where $a,b,u \in k$ are given by 
\begin{align}
a&=\frac{1}{8}\left(p^4+4p^2+2+p(p^2+2)\sqrt{p^2+4}\right),\nonumber\\
b&=\frac{5}{4}\left(p(p^2+2)(p^2+4)+(p^4+4p^2+2)\sqrt{p^2+4}\right)W_{p,r},
\label{abu}\\
u&=\frac{1}{16}\left((p^2+2)(p^4+4p^2+1)+p(p^2+1)(p^2+3)\sqrt{p^2+4}\right)
\nonumber
\end{align}
with $W_{p,r}=16(p^2+4)r^3 +4(p^2+4)r^2-4(19p+41)r-16p-199$ 
given as in (\ref{Wpr}). 

We obtain an isomorphism 
$f^{\ast}:\Eprast \rightarrow \Estast$ defined over $k$ 
such that the diagram 
\[
\begin{CD}
0 @>>> \Ker\nu^{\ast} @>>> \Eprast @>{\nu^{\ast}_{/\Q}}>> \Epr @>>> 0 \\
@.        @.						@V{f^{\ast}_{/k}}VV					@VV{f_{/k}}V @.\\
0 @>>> \Ker\phi^{\ast} @>>> \Estast @>>{\phi^{\ast}_{/k}}> \Est @>>>0\\
@.        @.						@V{g^{\ast}_{/F}}VV					@VV{g_{/F}}V @.\\
0 @>>> \Ker\lambda^{\ast} @>>> E_t^{\ast} @>>{\lambda^{\ast}_{/k}}> E_t @>>>0\\
\end{CD}
\]
commutes with exact rows. 
The $j$-invariants of $\Eprast$ and $\Estast$ are the same 
\begin{align*}
\frac{(p^2-12p+16)^3}{p-11}. 
\end{align*}
Therefore the isomorphism $f$ induces an injection 
\begin{align*}
\overline{f}:\Epr(k)/\nu^{\ast}(\Eprast(k)) 
\hookrightarrow \Est(k)/\phi^{\ast}(\Estast(k)).
\end{align*}

By \cite[Theorem 3.1]{KRS10} (see also Kida \cite[Remark 4.3]{Kid12}), 
there exists an injective homomorphism 
\begin{align*}
\overline{g}:\Est(k)/\phi^{\ast}(\Estast(k)) 
\hookrightarrow \Hom(\Gal(\overline{F}/F),\Ker\lambda^{\ast}(k)).
\end{align*}
Then we also obtain an injective homomorphism
\begin{align*}
\overline{g}\circ\overline{f}:\Epr(k)/\nu^{\ast}(\Eprast(k))\hookrightarrow
\Hom(\Gal(\overline{F}/F),\Ker\lambda^{\ast}(k)).
\end{align*}
Because the isogeny $\nu^{\ast}$ is defined over $\Q$, we get
\begin{align*}
\Epr(\Q)/\nu^{\ast}(\Eprast(\Q))\hookrightarrow\Hom(\Gal(\overline{F}/F),
\Ker\lambda^{\ast}(k)).
\end{align*}

Every point $P = (x(P),y(P)) \in \Epr(\Q)$ defines a Kummer extension 
\begin{align*}
L_P= F((\lambda^{\ast})^{-1}(g\circ f(P)))
\end{align*}
over $F$. 
In particular, via $(\ref{eqHM})$, we observe that 
\begin{align*}
L_P&=\mathrm{Spl}_k\left(\mathrm{Bru}
\left(\tfrac{1}{2}(p+\sqrt{p^2+4}),\tfrac{x(f(P))}{-4d};X\right)\right)\\
&=\mathrm{Spl}_\Q\left(\Lec(P;X)\right) 
\end{align*}
where $\displaystyle{\Lec(P;X)=\Lec\left(p,\frac{x(P)}{4(p^2+4)W};X\right)}$ 
as in Definition \ref{defLec}. 
Hence the group $\Epr(\Q)/\nu^{\ast}(\Eprast(\Q))$ classifies 
the isomorphism classes of $\mathrm{Spl}_\Q\left(\Lec(P;X)\right)$ 
with quartic subfield $F$ 
(see also \cite[Section 3]{KRS10}).\qed\\

By Theorem \ref{t1.4}, we have the following result 
by the multiplication-by-$2$ map of the elliptic curve $\Epr$: 

\begin{corollary}\label{c3.1}
For a $\Q$-rational point $P \in \Epr(\Q)$ and integer $n$ 
with $\gcd(n,5)=1$, 
${\rm Spl}_\Q(\Lec(P;X))={\rm Spl}_\Q(\Lec([n]P;X))$ 
where $\Lec(P;X)=\Lec \left(p,\frac{x(P)}{4(p^2+4)W};X\right)$ 
as in Definition \ref{defLec}. 
In particular, for $P=Q_0=(4r(p^2+4)W,2(p^2+4)W^2)$ and $n=2$, 
we have 
${\rm Spl}_\Q(\Lec(p,r;X))={\rm Spl}_\Q(\Lec(p,R;X))$ 
where 
\begin{align*}
R&=\frac{x([2]Q_0)}{4(p^2+4)W},\\
W&=16(p^2+4)r^3 +4(p^2+4)r^2-4(19p+41)r-16p-199,\\
x([2]Q_0)&=16(p^2+4)^2r^4+8(p^2+4)(19p+41)r^2+4(32p^3+398p^2+128p+1592)r\\
&\ \ \ +16p^3+560p^2+1622p+2477. 
\end{align*}
\end{corollary}

\begin{remark}
We can also verify 
${\rm Spl}_\Q(\Lec(p,r;X))={\rm Spl}_\Q(\Lec(p,R;X))$ 
in Corollary \ref{c3.1} 
by Hoshi and Miyake \cite{HM10} via multi-resolvent polynomials.
We take multi-resolvent polynomials $F^1_{a,a^\prime}$ and 
$F^2_{a,a^\prime}$ as in \cite[page 1071]{HM10} where 
$a=(s,t)$, $a^\prime=(s^\prime,t^\prime)$.
Using \cite[page 1078, Method 2]{HM10}, 
for
\begin{align*}
a&=\left( -\frac{1}{4}(5p+8r+2p^2r+(2pr+5)\sqrt{p^2+4}),
\frac{1}{2}(p+\sqrt{p^2+4})\right),\\
a^\prime&=\left(-\frac{1}{4}(5p+8R+2p^2R+(2pR+5)\sqrt{p^2+4}),
\frac{1}{2}(p+\sqrt{p^2+4})\right)
\end{align*}
via (\ref{eqHM}), we obtain that 
$\mathrm{Spl}_{\Q}(\Lec(p,r;X))=\mathrm{Spl}_{\Q}(\Lec(p,R;X))$ 
if and only if 
$F^1_{a,a^\prime}$ or $F^2_{a,a^\prime}$ has a linear factor over 
$k=\Q(\sqrt{p^2+4})$. 
Indeed, we can check that $F^2_{a,a^\prime}$ has a linear factor 
\begin{align*}
x+\frac{1+2r}{2}\sqrt{p^2+4}+\frac{p-1}{2}.
\end{align*}
\end{remark}

\section{Examples of Theorem \ref{t1.4}}

We will give two examples of Theorem \ref{t1.4} with 
$\mathcal{E}_{p,r}(\Q)/\nu^{\ast}(\mathcal{E}_{p,r}^{\ast}(\Q))
\simeq \Z/5\Z$ and $(\Z/5\Z)^{\oplus 2}$. 

\begin{example}[$p=1$ and $r=-3$ with $\mathcal{E}_{1,-3}(\Q)/\nu^{\ast}(\mathcal{E}_{1,-3}^{\ast}(\Q))\simeq \Z/5\Z$]\label{ex1}
We consider the case where $p=1$ and $r=-3$. 
The associated isogenous curves are 
\begin{align*}
\mathcal{E}_{1,-3} : y^2 &= x^3 - 7375x^2 - 2610750000x 
+ 68994507812500,\\
\mathcal{E}_{1,-3}^{\ast} : y^2 &= x^3 - 7375x^2 - 11313250000x 
- 5450566117187500
\end{align*}
with $j$-invariants 
$-\dfrac{5\cdot 29^3}{2^5}, -\dfrac{5^2}{2}$ respectively. 
Their Mordell--Weil groups are
\begin{align*}
\mathcal{E}_{1,-3}(\Q)&= \langle P_1,P_2 \rangle \simeq \Z^{\oplus 2},\\
\mathcal{E}_{1,-3}^{\ast}(\Q)&=\langle Q_1,Q_2\rangle \simeq \Z^{\oplus 2}
\end{align*}
where 
\begin{align*}
P_1&= (-53100, 6091750),& P_2&=(88500, 21756250),\\
Q_1&=(678500, 543906250),& Q_2&=(1452875, 1740500000).
\end{align*}
We can check $P_2=Q_0$ where $Q_0=(4r(p^2+4)W,2(p^2+4)W^2)$ 
which corresponds to $\Lec(1,-3;X)$. 
The isogeny $\nu^{\ast}:\mathcal{E}_{1,-3}^{\ast} \rightarrow 
\mathcal{E}_{1,-3}$ is given by 
\begin{align*}
\nu^{\ast}(Q_1)&=P_1-2P_2,\\
\nu^{\ast}(Q_2)&=-P_1-3P_2.
\end{align*}
Hence the image of $\nu^{\ast}$ is given by 
\begin{align*}
\nu^{\ast}(\mathcal{E}_{1,-3}^{\ast})=\langle P_1-2P_2, 5P_2\rangle.
\end{align*}
We conclude that 
$\mathcal{E}_{1,-3}(\Q)/\nu^{\ast}(\mathcal{E}_{1,-3}^{\ast}(\Q))
=\langle\overline{P_2}\rangle\simeq \Z/5\Z$. 
Thus there exists exactly one isomorphism class 
of Lecachux's polynomials. 
We have 
\begin{align*}
\Spl{\Lec(1,3:X)}=\Spl{\Lec([n]P_2:X)}
=\Spl{\Lec\left(1,\frac{x([n]P_2)}{4(p^2+4)W};X\right)}
\end{align*}
where $\gcd(n,5)=1$. 
For example, for $n=1$, $2$, $3$, $4$, we have 
\begin{align*}
\frac{x([n]P_2)}{4(p^2+4)W}=
-3, 
\frac{-263}{236}, 
\frac{4849}{39605}, 
\frac{2034016227}{1036798976}
\end{align*}
respectively. 
We can check this example by Sage (\cite{Sage}) as follows:\\

{\small 
\begin{verbatim}
sage: p=1;r=-3;e=p^2+4;W=16*e*r^3+4*e*r^2-4*(19*p+41)*r-16*p-199;
sage: Epr=EllipticCurve([0,e*W,0,-4*(19*p+41)*e*W^2,-4*e^2*(16*p+199)*W^3]);Epr
Elliptic Curve defined by y^2 = x^3 - 7375*x^2 - 2610750000*x 
+ 68994507812500 over Rational Field
sage: factor(Epr.j_invariant())
-1 * 2^-5 * 5 * 29^3
sage: Epr.torsion_points()
[(0 : 1 : 0)]
sage: EprGen=Epr.gens();EprGen
[(-53100 : 6091750 : 1), (88500 : 21756250 : 1)]
sage: f=list(factor(Epr.division_polynomial(5)))[0][0];f
x - 44250
sage: x=44250
sage: u=x^3 - 7375*x^2 - 2610750000*x + 68994507812500;u
25672375000000
sage: factor(u)
2^6 * 5^9 * 59^3
sage: R.<x>=PolynomialRing(QQ)
sage: K.<a>=NumberField(x^2-u);K
Number Field in a with defining polynomial x^2 - 25672375000000
sage: EprK=Epr.base_extend(K);EprK
Elliptic Curve defined by y^2 = x^3 + (-7375)*x^2 + (-2610750000)*x 
+ 68994507812500 over Number Field in a with defining polynomial x^2 
- 25672375000000
sage: P=EprK(44250,a);P
(44250 : a : 1)
sage: phi=EprK.isogeny(P)
sage: EprKast=phi.codomain()
sage: Eprast=EprKast.base_extend(QQ);Eprast
Elliptic Curve defined by y^2 = x^3 - 7375*x^2 - 11313250000*x 
- 5450566117187500 over Rational Field
sage: Eprast.j_invariant()
-25/2
sage: (p^2-12*p+16)^3/(p-11)
-25/2
sage: Eprast.torsion_points()
[(0 : 1 : 0)]
sage: EprastGen=Eprast.gens();EprastGen
[(678500 : 543906250 : 1), (1452875 : 1740500000 : 1)]
sage: P1=EprGen[0];P1
(-53100 : 6091750 : 1)
sage: P2=EprGen[1];P2
(88500 : 21756250 : 1)
sage: Q1=EprastGen[0];Q1
(678500 : 543906250 : 1)
sage: Q2=EprastGen[1];Q2
(1452875 : 1740500000 : 1)
sage: phiast=phi.dual()
sage: Q1_img=Epr(phiast(Q1))
sage: Q2_img=Epr(phiast(Q2))
sage: for i in range(-5,6):
....:     for j in range(-5,6):
....:         if i*P1+j*P2==Q1_img:
....:             [i,j]
....:
[1, -2]
sage:  for i in range(-5,6):
....:     for j in range(-5,6):
....:         if i*P1+j*P2==Q2_img:
....:             [i,j]
....:
[-1, -3]
sage: Q0=Epr(4*r*(p^2+4)*W,2*(p^2+4)*W^2);Q0
(88500 : 21756250 : 1)
sage: P2==Q0
True
sage: for i in range(1,5):
....:     (i*P2)[0]/(4*(p^2+4)*W)
....:
-3
-263/236
4849/39605
2034016227/1036798976
\end{verbatim}
}
\end{example}

\begin{example}[$p=2$ and $r=-15$ with $\mathcal{E}_{2,-15}(\Q)/\nu^{\ast}(\mathcal{E}_{2,-15}^{\ast}(\Q))\simeq (\Z/5\Z)^{\oplus 2}$]\label{ex2}
We consider the case where $p=2$, $r=-15$. 
The associated isogenous curves are
\begin{align*}
\mathcal{E}_{2,-15} :y^2 &= 
x^3 - 3362328x^2 - 446557358393568x + 4390381057572915584256,\\
\mathcal{E}_{2,-15}^{\ast} : y^2 &= 
x^3 - 3362328x^2 + 1181398581066528x - 243295532112514685688576
\end{align*}
with $j$-invariants 
$-\dfrac{2^6\cdot 239^3}{3^{10}}, \dfrac{2^6}{3^2}$ 
respectively. 
Their Mordell--Weil groups are
\begin{align*}
\mathcal{E}_{2,-15}(\Q)&=
\langle P_{\Tor}\rangle\oplus\langle P_1,P_2,P_3 \rangle 
\simeq \Z/2\Z \oplus \Z^{\oplus 3},\\
\mathcal{E}_{2,-15}^{\ast}(\Q)&=
\langle Q_{\Tor}\rangle\oplus\langle Q_1,Q_2,Q_3 \rangle 
\simeq \Z/2\Z \oplus \Z^{\oplus 3}
\end{align*}
where 
\begin{align*}
P_{\mathrm{tor}}&=(-23536296,0),& 
P_1&=\left( \frac{1213850592}{121}, \frac{32104365187824}{1331} \right),\\
P_2&=(12954852, 14669441496),& 
P_3&=(24185016, 75959770464),\\
Q_{\mathrm{tor}}&=(57159576,0),& 
Q_1&=\left( \frac{9662338144}{169} , \frac{26786536642000}{2197} \right),\\
Q_2&=(58184676 , 105083001000),& 
Q_3&=\left( \frac{15400097496}{121} , \frac{1841522732064000}{1331} \right).
\end{align*}
The isogeny $\nu^{\ast}:\mathcal{E}_{2,-15}^{\ast} \rightarrow 
\mathcal{E}_{2,-15}$ is given by 
\begin{align*}
\nu^{\ast}(Q_{\Tor})&=P_{\Tor},\\
\nu^{\ast}(Q_1)&=-P_1+2P_2+2P_3,\\
\nu^{\ast}(Q_2)&=P_{\Tor}-2P_1-P_2-P_3,\\
\nu^{\ast}(Q_3)&=-2P_1+4P_2-P_3.
\end{align*}
Hence we obtain the image 
\begin{align*}
\nu^{\ast}(\mathcal{E}_{2,-15}^{\ast})
=\langle P_{\Tor},P_1+2P_2+2P_3, 5P_2, 5P_3\rangle
\end{align*}
and conclude that 
$\mathcal{E}_{2,-15}(\Q)/\nu^{\ast}(\mathcal{E}_{2,-15}^{\ast}(\Q))
=\langle\overline{P_2},\overline{P_3}\rangle
\simeq (\Z/5\Z)^{\oplus 2}$. 
There exist $6$ subgroups of order $5$ 
in $\mathcal{E}_{2,-15}(\Q)/\nu^{\ast}(\mathcal{E}_{2,-15}^{\ast}(\Q))
\simeq (\Z/5\Z)^{\oplus 2}$ 
which correspond to the $6$ isomorphism classes 
\begin{align*}
\Lec(P_2-2P_3;X)&=\Lec\left(p,\frac{x(P_2-2P_3)}{4(p^2+4)W};X\right)
=\Lec\left(2,-\frac{6826408529368884683}{114084259282587016};X\right),\\
\Lec(P_2-P_3;X)&=\Lec\left(p,\frac{x(P_2-P_3)}{4(p^2+4)W};X\right)
=\Lec\left(2,-\frac{5293745}{2271049};X\right),\\
\Lec(P_2;X)&=\Lec\left(p,\frac{x(P_2)}{4(p^2+4)W};X\right)
=\Lec\left(2,-\frac{131}{136};X\right),\\
\Lec(P_2+P_3;X)&=\Lec\left(p,\frac{x(P_2+P_3)}{4(p^2+4)W};X\right)
=\Lec\left(2,\frac{157}{529};X\right),\\
\Lec(P_2+2P_3;X)&=\Lec\left(p,\frac{x(P_2+2P_3)}{4(p^2+4)W};X\right)
=\Lec\left(2,\frac{9701177386741}{7753965979144};X\right),\\
\Lec(P_3;X)&=\Lec\left(p,\frac{x(P_3)}{4(p^2+4)W};X\right)
=\Lec\left(2,-\frac{19759}{10988};X\right). 
\end{align*}
with the quartic subfield 
\begin{align*}
F=\Q\left(\sqrt{-233495-\frac{326893}{2}\sqrt{2}}\right). 
\end{align*}
Since ${\rm Lec}(2,-15;X)$ corresponds to the rational point 
\begin{align*}
Q_0=(4r(p^2+4)W,2(p^2+4)W^2)=(201739680,2826312394896)
=P_{\Tor}-P_1-P_3
\end{align*}
and $\langle\overline{Q_0}\rangle=\langle\overline{P_2-2P_3}\rangle$ in 
$\mathcal{E}_{2,-15}(\Q)/\nu^{\ast}(\mathcal{E}_{2,-15}^{\ast}(\Q))$, 
we have 
\begin{align*}
\Spl{\Lec(Q_0;X)}
=\Spl{\Lec(2,-15;X)}=\Spl{\Lec(P_2-2P_3;X)}.
\end{align*}
We can check this example by Sage (\cite{Sage}) as follows:\\

{\small 
\begin{verbatim}
sage: p=2;r=-15;e=p^2+4;W=16*e*r^3+4*e*r^2-4*(19*p+41)*r-16*p-199;
sage: Epr=EllipticCurve([0,e*W,0,-4*(19*p+41)*e*W^2,-4*e^2*(16*p+199)*W^3]);Epr
Elliptic Curve defined by y^2 = x^3 - 3362328*x^2 - 446557358393568*x 
+ 4390381057572915584256 over Rational Field
sage: factor(Epr.j_invariant())
-1 * 2^6 * 3^-10 * 239^3
sage: Epr.torsion_points()
[(-23536296 : 0 : 1), (0 : 1 : 0)]
sage: EprGen=Epr.gens();EprGen
[(1213850592/121 : 32104365187824/1331 : 1),
 (12954852 : 14669441496 : 1),
 (24185016 : 75959770464 : 1)]
sage: f=list(factor(Epr.division_polynomial(5)))[0][0];f
x^2 - 13449312*x - 231757616381472/5
sage: R.<x>=PolynomialRing(QQ)
sage: F.<a>=NumberField(x^2-13449312*x-231757616381472/5);F
Number Field in a with defining polynomial x^2 - 13449312*x - 231757616381472/5
sage: x=a
sage: u=x^3 - 3362328*x^2 - 446557358393568*x + 4390381057572915584256;u
-1322714200811328/5*a + 24289640656182623881728/5
sage: R2.<y>=PolynomialRing(F)
sage: K.<b>=F.extension(y^2-u);K
Number Field in b with defining polynomial y^2 + 1322714200811328/5*a 
- 24289640656182623881728/5 over its base field
sage: EprK=Epr.base_extend(K);EprK
Elliptic Curve defined by y^2 = x^3 + (-3362328)*x^2 + (-446557358393568)*x 
+ 4390381057572915584256 over Number Field in b with defining polynomial 
y^2 + 1322714200811328/5*a - 24289640656182623881728/5 over its base field
sage: P=EprK(a ,b);P
(a : b : 1)
sage: phi=EprK.isogeny(P)
sage: EprKast=phi.codomain()
sage: Eprast=EprKast.base_extend(QQ);Eprast
Elliptic Curve defined by y^2 = x^3 - 3362328*x^2 + 1181398581066528*x 
- 243295532112514685688576 over Rational Field
sage: Eprast.j_invariant()
64/9
sage: (p^2-12*p+16)^3/(p-11)
64/9
sage: Eprast.torsion_points()
[(0 : 1 : 0), (57159576 : 0 : 1)]
sage: EprastGen=Eprast.gens();EprastGen
[(9662338144/169 : 26786536642000/2197 : 1),
 (58184676 : 105083001000 : 1),
 (15400097496/121 : 1841522732064000/1331 : 1)]
sage: Ptor=Epr.torsion_points()[0];Ptor
(-23536296 : 0 : 1)
sage: P1=EprGen[0];P1
(1213850592/121 : 32104365187824/1331 : 1)
sage: P2=EprGen[1];P2
(12954852 : 14669441496 : 1)
sage: P3=EprGen[2];P3
(24185016 : 75959770464 : 1)
sage: Qtor=Eprast.torsion_points()[1];Qtor
(57159576 : 0 : 1)
sage: Q1=EprastGen[0];Q1
(9662338144/169 : 26786536642000/2197 : 1)
sage: Q2=EprastGen[1];Q2
(58184676 : 105083001000 : 1)
sage: Q3=EprastGen[2];Q3
(15400097496/121 : 1841522732064000/1331 : 1)
sage: phiast=phi.dual()
sage: Qtor_img=Epr(phiast(Qtor))
sage: Q1_img=Epr(phiast(Q1))
sage: Q2_img=Epr(phiast(Q2))
sage: Q3_img=Epr(phiast(Q3))
sage: Qtor_img==Ptor
True
sage: for i in range(0,2):
....:     for j in range(-5,6):
....:         for k in range(-5,6):
....:             for l in range(-5,6):
....:                 if Q1_img==i*Ptor+j*P1+k*P2+l*P3:
....:                     [i,j,k,l]
....:
[0, -1, 2, 2]
sage: for i in range(0,2):
....:     for j in range(-5,6):
....:         for k in range(-5,6):
....:             for l in range(-5,6):
....:                 if Q2_img==i*Ptor+j*P1+k*P2+l*P3:
....:                     [i,j,k,l]
....:
[1, -2, -1, -1]
sage: for i in range(0,2):
....:     for j in range(-5,6):
....:         for k in range(-5,6):
....:             for l in range(-5,6):
....:                 if Q3_img==i*Ptor+j*P1+k*P2+l*P3:
....:                     [i,j,k,l]
....:
[0, -2, 4, -1]
sage: Q0=Epr(4*r*(p^2+4)*W,2*(p^2+4)*W^2);Q0
(201739680 : 2826312394896 : 1)
sage: for i in range(0,2):
....:     for j in range(-5,6):
....:         for k in range(-5,6):
....:             for l in range(-5,6):
....:                 if Q0==i*Ptor+j*P1+k*P2+l*P3:
....:                     [i,j,k,l]
....:
[1, -1, 0, -1]
sage: for i in range(-2,3):
....:     (P2+i*P3)[0]/(4*(p^2+4)*W)
....:
-6826408529368884683/114084259282587016
-5293745/2271049
-131/136
157/529
9701177386741/7753965979144
sage: P3[0]/(4*(p^2+4)*W)
-19759/10988
\end{verbatim}
}
\end{example}

\section{{Degenerate cases with $p=0$ in Theorem \ref{t1.4}}}
By Hoshi and Miyake \cite[Lemma 7.3]{HM10}, 
for $(p,r) \in \Q^2$, 
there exists $b \in \Q$ such that $b^2=p^2+4$ 
if and only if $G_{p,r} \le D_5$ 
where $G_{p,r}=\Gal(\Lec(p,r;X)/\Q)$. 
Moreover, in this case where $G_{p,r} \le D_5$, 
by (\ref{eqHM}), 
the splitting fields of $\Lec(p,r;X)$ and of $\Bru(t,s;X)$ over $\Q$ 
coincide 
where 
\begin{align*}
s=-\frac{5p+8r+2p^2r+(2pr+5)b}{4}, t=\frac{p+b}{2}
\end{align*} 
(see also Kida, Rikuna and Sato \cite[Example 4.1]{KRS10}). 

We will give two examples of the case $p=0$ $(b^2=p^2+4=4)$ 
with $G_{p,r}\simeq D_5$ and $C_5$. 

\begin{example}[$p=0, r=-\frac{5}{4}$ with $G_{p,r}\simeq D_5$]\label{ex3}

We consider the case where $p=0$ and $r=-\frac{5}{4}$. 
Then we see $\Lec\left(0,-\frac{5}{4};X\right)=\Bru(1,0;X)$ 
and $G_{p,r}\simeq D_5$. 
The associated isogenous curves are 
\begin{align*}
\mathcal{E}_{0,-\frac{5}{4}}:y^2&=x^3 - 376x^2 - 5796416x + 10578317824,\\
\mathcal{E}_{0,-\frac{5}{4}}^\ast:y^2&=x^3-376x^2-117766208x-1154206105088
\end{align*}
with $j$-invariants 
$-\dfrac{2^{12}\cdot 31^3}{11^5}, -\dfrac{2^{12}}{11}$ 
respectively. 
Their Mordell--Weil groups are
\begin{align*}
\mathcal{E}_{0,-\frac{5}{4}}(\Q)&=\langle P_1, P_2 \rangle
\simeq \Z^{\oplus 2},\\
\mathcal{E}_{0,-\frac{5}{4}}^\ast(\Q)&=\langle Q_1, Q_2 \rangle 
\simeq \Z^{\oplus 2}
\end{align*}
where 
\begin{align*}
P_1&=(-2632,70688),& &P_2=(1880,70688),\\
Q_1&=(20492,2209000),& &Q_2=(43992,8836000). 
\end{align*}
We can check $P_2=Q_0$ where $Q_0=(4r(p^2+4)W,2(p^2+4)W^2)$ 
which corresponds to $\Lec(0,-\frac{5}{4};X)$. 
The isogeny 
$\nu^{\ast}:\mathcal{E}_{0,-\frac{5}{4}}^\ast \rightarrow 
\mathcal{E}_{0,-\frac{5}{4}}$ 
is given by 
\begin{align*}
\nu^{\ast}(Q_1)&=-2P_1-3P_2,\\ 
\nu^{\ast}(Q_2)&=P_1-P_2.
\end{align*}
Hence we have 
\begin{align*}
\nu^{\ast}(\mathcal{E}_{0,-\frac{5}{4}}^\ast)=\langle P_1-P_2, 5P_2 \rangle.
\end{align*}
Then we conclude that 
$\mathcal{E}_{0,-\frac{5}{4}}(\Q)
/\nu^{\ast}(\mathcal{E}_{0,-\frac{5}{4}}^{\ast}(\Q))
=\langle\overline{P_2}\rangle\simeq \Z/5\Z$. 

On the other hand, we have 
\begin{align*}
E_{1,0}:y^2=x^3 + 1316x^2 + 212064x + 78074896
\end{align*}
with the Mordell--Weil group 
\begin{align*}
E_{1,0}(\Q)=\langle P_1^\prime,\ P_2^\prime \rangle \simeq \Z^{\oplus2}
\end{align*}
where 
\begin{align*}
P_1^\prime=(-1128,8836), P_2^\prime=(0,8836)
\end{align*}
(see also Kida, Rikuna and Sato \cite[Example 4.1]{KRS10}). 
We see $P_2^\prime=P_0$
where 
$P_0=(-4ds,4d^2)$ 
which corresponds to $\Bru(1,0;X)$. 
The $j$-invariant of $E_{1,0}$ and of $\mathcal{E}_{0,-\frac{5}{4}}$ are 
the same $-\dfrac{2^{12}\cdot 31^3}{11^5}$. 
Indeed, by (\ref{abu}) in the proof of Theorem \ref{t1.4}, 
we get the isomorphism 
\begin{align*}
f:&\ \mathcal{E}_{0,-\frac{5}{4}}\rightarrow E_{1,0},\\
&\ (x,y) \mapsto \left(\frac{1}{4}x-470,\ \frac{1}{8}y\right)
\end{align*}
with 
\begin{align*}
f(P_1)=P_1^\prime\ \text{and}\ f(P_2)=P_2^\prime.
\end{align*}
Hence it follows from Theorem \ref{t1.4} that 
\begin{align*}
\Spl{\Lec(0,-\tfrac{5}{4};X)}=
\Spl{\Lec(P_2;X)}=\Spl{\Bru(P_2^\prime;X)}=\Spl{\Bru(1,0;X)}. 
\end{align*} 
We can check this example by Sage (\cite{Sage}) as follows:\\

{\small 
\begin{verbatim}
sage: p=0;r=-5/4;e=p^2+4;W=16*e*r^3+4*e*r^2-4*(19*p+41)*r-16*p-199;
sage: Epr=EllipticCurve([0,e*W,0,-4*(19*p+41)*e*W^2,-4*e^2*(16*p+199)*W^3]);Epr
Elliptic Curve defined by y^2 = x^3 - 376*x^2 - 5796416*x + 10578317824 
over Rational Field
sage: factor(Epr.j_invariant())
-1 * 2^12 * 11^-5 * 31^3
sage: Epr.torsion_points()
[(0 : 1 : 0)]
sage: EprGen=Epr.gens();EprGen
[(-2632 : 70688 : 1), (1880 : 70688 : 1)]
sage: factor(Epr.division_polynomial(5))
(5) * (x + 3384) * (x + 11656) * (x^2 - 1504*x - 14278976/5) 
* (x^4 - 12784*x^3 + 81432576*x^2 - 141079675904*x + 91761112969216) 
* (x^4 - 2256*x^3 + 7917056*x^2 - 22645042176*x + 39794462191616)
sage: R.<x>=PolynomialRing(QQ)
sage: F.<a>=NumberField(x^2-1504*x-14278976/5);F
Number Field in a with defining polynomial x^2 - 1504*x - 14278976/5
sage: x=a
sage: u=x^3 - 376*x^2 - 5796416*x + 10578317824;u
-6220544/5*a + 68998274048/5
sage: R2.<y>=PolynomialRing(F)
sage: K.<b>=F.extension(y^2-u);K
Number Field in b with defining polynomial y^2 + 6220544/5*a - 68998274048/5 
over its base field
sage: EprK=Epr.base_extend(K);EprK
Elliptic Curve defined by y^2 = x^3 + (-376)*x^2 + (-5796416)*x + 10578317824 
over Number Field in b with defining polynomial y^2 + 6220544/5*a 
- 68998274048/5 over its base field
sage: P=EprK(a,b);P
(a : b : 1)
sage: phi=EprK.isogeny(P)
sage: EprKast=phi.codomain()
sage: Eprast=EprKast.base_extend(QQ);Eprast
Elliptic Curve defined by y^2 = x^3 - 376*x^2 - 117766208*x - 1154206105088 
over Rational Field
sage: factor(Eprast.j_invariant())
-1 * 2^12 * 11^-1
sage: factor((p^2-12*p+16)^3/(p-11))
-1 * 2^12 * 11^-1
sage: Eprast.torsion_points()
[(0 : 1 : 0)]
sage: EprastGen=Eprast.gens();EprastGen
[(20492 : 2209000 : 1), (43992 : 8836000 : 1)]
sage: P1=EprGen[0];P1
(-2632 : 70688 : 1)
sage: P2=EprGen[1];P2
(1880 : 70688 : 1)
sage: Q1=EprastGen[0];Q1
(20492 : 2209000 : 1)
sage: Q2=EprastGen[1];Q2
(43992 : 8836000 : 1)
sage: phiast=phi.dual()
sage: Q1_img=Epr(phiast(Q1))
sage: Q2_img=Epr(phiast(Q2))
sage: for i in range(-5,6):
....:     for j in range(-5,6):
....:         if i*P1+j*P2==Q1_img:
....:             [i,j]
....:
[-2, -3]
sage: for i in range(-5,6):
....:     for j in range(-5,6):
....:         if i*P1+j*P2==Q2_img:
....:             [i,j]
....:
[1, -1]
sage: Q0=Epr(4*r*(p^2+4)*W,2*(p^2+4)*W^2);Q0
(1880 : 70688 : 1)
sage: Q0==P2
True
sage: t=1;s=0;d=-4*s^3+(t^2-30*t+1)*s^2+2*t*(3*t+1)*(4*t-7)*s-t*(4*t^4-4*t^3-40*
....: t^2+91*t-4);d
-47
sage: Ets=EllipticCurve([0,d*(t^2-30*t+1),0,-8*d^2*t*(3*t+1)*(4*t-7),-16*d^3*(4*
....: t^4-4*t^3-40*t^2+91*t-4)]);Ets
Elliptic Curve defined by y^2 = x^3 + 1316*x^2 + 212064*x + 78074896 
over Rational Field
sage: Ets.torsion_points()
[(0 : 1 : 0)]
sage: EtsGen=Ets.gens();EtsGen
[(-1128 : 8836 : 1), (0 : 8836 : 1)]
sage: P1p=EtsGen[0];P1p
(-1128 : 8836 : 1)
sage: P2p=EtsGen[1];P2p
(0 : 8836 : 1)
sage: P0=Ets(-4*d*s,4*d^2);P0
(0 : 8836 : 1)
sage: P0==P2p
True
sage: fP1=Ets((1/4)*P1[0]-470,P1[1]/8);fP1
(-1128 : 8836 : 1)
sage: fP2=Ets((1/4)*P2[0]-470,P2[1]/8);fP2
(0 : 8836 : 1)
sage: fP1==P1p
True
sage: fP2==P2p
True
\end{verbatim}
}
\end{example}

\begin{example}[$p=0, r=\frac{31}{4}$ with $G_{p,r}\simeq C_5$]\label{ex4}

We treat the degenerate case with 
$\Gal(\Lec(p,r;X)/\Q)=G_{p,r}\simeq C_5$, i.e. 
$\Q(p,r)(\sqrt{D_{p,r}})=\Q$. 
Take $p=0$, $r=\frac{31}{4}$ with 
$\Lec\left(0,\frac{31}{4};X\right)=\Bru(1,-18;X)$ 
and $G_{p,r}\simeq C_5$. 
The associated isogenous curves are
\begin{align*}
\mathcal{E}_{0,\frac{31}{4}}:  y^2 = x^3 + 117128x^2 - 562477703744x 
- 319768190447349248,
\end{align*}
\begin{align*}
\mathcal{E}_{0,\frac{31}{4}}^{\ast} : y^2 = x^3 + 117128x^2 
- 11427900663872x + 348900840159964537
\end{align*}
with $j$-invariants 
$-\dfrac{2^{12}\cdot 31^3}{11^5}, -\dfrac{2^{12}}{11}$ 
respectively. 
Their Mordell--Weil groups are
\begin{align*}
\mathcal{E}_{0,\frac{31}{4}}(\Q)&= \langle P_{\Tor} \rangle \simeq \Z/5\Z,\\ 
\mathcal{E}_{0,\frac{31}{4}}^{\ast}(\Q)&=\langle Q_{\Tor}\rangle\simeq \Z/5\Z
\end{align*}
where
\begin{align*}
P_{\Tor}=(3630968,6859484192), 
Q_{\Tor}=(-1991176,-7086244000).
\end{align*}
We see that $P_{\Tor}=Q_0$ where $Q_0=(4r(p^2+4)W,2(p^2+4)W^2)$ 
which corresponds to $\Lec(0,\frac{31}{4};X)$. 
The isogeny $\nu^{\ast}: \mathcal{E}_{0,\frac{31}{4}}^{\ast}
\rightarrow \mathcal{E}_{0,\frac{31}{4}}$ is given by 
\begin{align*}
\nu^{\ast}(Q_{\Tor})=O
\end{align*}
and hence we conclude that 
$\mathcal{E}_{0,\frac{31}{4}}(\Q)/\nu^{\ast}
(\mathcal{E}_{0,\frac{31}{4}}^{\ast}(\Q))=\langle P_{\Tor}
\rangle\simeq \Z/5\Z$. 

On the other hand, we have 
\begin{align*}
E_{1,-18}:y^2=x^3-409948x^2 + 20578452576x -2360098139294192
\end{align*}
with the Mordell--Weil group 
\begin{align*}
E_{1,-18}(\Q)=\langle P_{\Tor}^\prime\rangle \simeq\Z/5\Z
\end{align*}
where 
\begin{align*}
P_{\Tor}^\prime=(1054152,857435524)
\end{align*}
(see also Kida, Rikuna and Sato \cite[Example 4.3]{KRS10}). 
We also have $P_{\Tor}^\prime=P_0$
where 
$P_0=(-4ds,4d^2)$ which corresponds to $\Bru(1,-18;X)$. 
The $j$-invariant of $E_{1,-18}$ and of $\mathcal{E}_{0,\frac{31}{4}}$ 
are the same $-\dfrac{2^{12}\cdot 31^3}{11^5}$. 
Indeed, by (\ref{abu}) in the proof of Theorem \ref{t1.4}, 
we obtain the isomorphism 
\begin{align*}
f:&\ \mathcal{E}_{0,\frac{31}{4}}\rightarrow E_{1,-18},\\
&\ (x,y) \mapsto \left(\frac{1}{4}x+146410,\ \frac{1}{8}y\right)
\end{align*}
with 
\begin{align*}
f(P_{\Tor})=P_{\Tor}^\prime. 
\end{align*}
Hence it follows from Theorem \ref{t1.4} that 
\begin{align*}
\Spl{\Lec(0,\tfrac{31}{4};X)}=
\Spl{\Lec(P_{\Tor};X)}=\Spl{\Bru(P_{\Tor}^\prime;X)}=\Spl{\Bru(1,-18;X)}. 
\end{align*} 
We can check this example by Sage (\cite{Sage}) as follows:\\

{\small 
\begin{verbatim}
sage: p=0;r=31/4;e=p^2+4;W=16*e*r^3+4*e*r^2-4*(19*p+41)*r-16*p-199;
sage: Epr=EllipticCurve([0,e*W,0,-4*(19*p+41)*e*W^2,-4*e^2*(16*p+199)*W^3]);Epr
Elliptic Curve defined by y^2 = x^3 + 117128*x^2 - 562477703744*x 
- 319768190447349248 ver Rational Field
sage: factor(Epr.j_invariant())
-1 * 2^12 * 11^-5 * 31^3
sage: Epr.gens()
[]
sage: Epr.torsion_points()
[(0 : 1 : 0),
 (1054152 : -623589472 : 1),
 (1054152 : 623589472 : 1),
 (3630968 : -6859484192 : 1),
 (3630968 : 6859484192 : 1)]
sage: factor(Epr.division_polynomial(5))
(5) * (x - 3630968) * (x - 1054152) * (x^2 + 468512*x - 1385615806784/5) 
* (x^4 + 70276*x^3 + 768262229504*x^2 + 684528890103370752*x 
+ 374726296200692418768896) * (x^4 + 398352*x^3 + 7902125789184*x^2 
+ 4264647122850577408*x + 864072539355790504554496)
sage: R.<x>=PolynomialRing(QQ)
sage: F.<a>=NumberField(x^2+468512*x-1385615806784/5);F
Number Field in a with defining polynomial x^2 + 468512*x - 1385615806784/5
sage: x=a
sage: u=x^3+117128*x^2-562477703744*x-319768190447349248;u
-603634608896/5*a - 2085724176887735296/5
sage: R2.<y>=PolynomialRing(F)
sage: K.<b>=F.extension(y^2-u);K
Number Field in b with defining polynomial y^2 + 603634608896/5*a 
+ 20857241768877352965 over its base field
sage: EprK=Epr.base_extend(K);EprK
Elliptic Curve defined by y^2 = x^3 + 117128*x^2 + (-562477703744)*x 
+ (-319768190447349248) over Number Field in b with defining polynomial 
y^2 + 603634608896/5*a + 2085724176887735296/5 over its base field
sage: P=EprK(a,b);P
(a : b : 1)
sage: phi=EprK.isogeny(P)
sage: EprKast=phi.codomain()
sage: Eprast=EprKast.base_extend(QQ);Eprast
Elliptic Curve defined by y^2 = x^3 + 117128*x^2 - 11427900663872*x 
+ 3489008401599645376 over Rational Field
sage: factor(Eprast.j_invariant())
-1 * 2^12 * 11^-1
sage: factor((p^2-12*p+16)^3/(p-11))
-1 * 2^12 * 11^-1
sage: EprastGen=Eprast.gens();EprastGen
[]
sage: Eprast.torsion_points()
[(-1991176 : -7086244000 : 1),
 (-1991176 : 7086244000 : 1),
 (0 : 1 : 0),
 (3865224 : -7086244000 : 1),
 (3865224 : 7086244000 : 1)]
sage: Qtor1=Eprast.torsion_points()[0];Qtor1
(-1991176 : -7086244000 : 1)
sage: phiast=phi.dual()
sage: Qtor1_img=Epr(phiast(Qtor1));Qtor1_img
(0 : 1 : 0)
sage: Q0=Epr(4*r*(p^2+4)*W,2*(p^2+4)*W^2);Q0
(3630968 : 6859484192 : 1)
sage: Ptor4=Epr.torsion_points()[4];Ptor4
(3630968 : 6859484192 : 1)
sage: Q0==Ptor4
True
sage: t=1;s=-18;d=-4*s^3+(t^2-30*t+1)*s^2+2*t*(3*t+1)*(4*t-7)*s-t*(4*t^4-4*t^3-4
....: 0*t^2+91*t-4);d
14641
sage: Ets=EllipticCurve([0,d*(t^2-30*t+1),0,-8*d^2*t*(3*t+1)*(4*t-7),-16*d^3*(4*
....: t^4-4*t^3-40*t^2+91*t-4)]);Ets
Elliptic Curve defined by y^2 = x^3 - 409948*x^2 + 20578452576*x 
- 2360098139294192 over Rational Field
sage: Ets.gens()
[]
sage: Ets.torsion_points()
[(0 : 1 : 0),
 (409948 : -77948684 : 1),
 (409948 : 77948684 : 1),
 (1054152 : -857435524 : 1),
 (1054152 : 857435524 : 1)]
sage: Ptor4p=Ets.torsion_points()[4];Ptor4p
(1054152 : 857435524 : 1)
sage: P0=Ets(-4*d*s,4*d^2);P0
(1054152 : 857435524 : 1)
sage: P0==Ptor4p
True
sage: fPtor4=Ets((1/4)*Ptor4[0]+146410,Ptor4[1]/8);fPtor4
(1054152 : 857435524 : 1)
sage: fPtor4==Ptor4p
True
\end{verbatim}
}
\end{example}\vspace*{1mm}
\begin{acknowledgment}
The authors thank the referee for helpful comments. 
\end{acknowledgment}

\end{document}